\documentclass{ifacconf}

\usepackage{amsmath}
\usepackage{amssymb}

\usepackage{graphicx}
\usepackage{float}
\PassOptionsToPackage{hyphens}{url}\usepackage{url}
\usepackage{algorithm}
\usepackage{algpseudocode}

\usepackage{stackengine}
\setstackEOL{\\}
\makeatletter
\def\BState{\State\hskip-\ALG@thistlm}
\makeatother

\usepackage{flushend}

\makeatletter
\def\endfigure{\end@float}
\def\endtable{\end@float}
\makeatother

\let\ifacconfcaptionwidth\captionwidth
\usepackage[caption=false]{subfig}
\let\captionwidth\ifacconfcaptionwidth
\usepackage{stfloats}
\newtheorem{assumption}{Assumption}
\newtheorem{definition}{Definition}

\newcommand{\mbb}{\mathbb}
\renewcommand{\Re}{\mbb{R}}

\newcommand{\set}[2]{\left\{ #1\ \left| \ #2 \right. \right\}}

\newcommand{\innerprod}[2]{\langle{#1},{#2}\rangle}

\usepackage{natbib}        % required for bibliography

\begin{document}
\begin{frontmatter}

\title{Gradient-Bounded Dynamic Programming with Submodular and Concave Extensible Value Functions\thanksref{footnoteinfo}} 

\thanks[footnoteinfo]{Research is supported by SIA Food Union Management.}

\author[First]{Denis Lebedev, Paul Goulart \& Kostas Margellos} 
%\author[Second]{Paul Goulart} 
%\author[Third]{Kostas Margellos}

\address[First]{Department of Engineering Science, University of Oxford, Oxford, OX1 3PJ, United Kingdom (e-mail: \{denis.lebedev, paul.goulart, kostas.margellos\}@eng.ox.ac.uk)}
%\address[Second]{Colorado State University, 
%   Fort Collins, CO 80523 USA (e-mail: author@lamar. colostate.edu)}
%\address[Third]{Electrical Engineering Department, 
%   Seoul National University, Seoul, Korea, (e-mail: author@snu.ac.kr)}

\begin{abstract}                % Abstract of not more than 250 words.
We consider dynamic programming problems with finite, discrete-time horizons and prohibitively high-dimensional, discrete state-spaces for direct computation of the value function from the Bellman equation. For the case that the value function of the dynamic program is concave extensible and submodular in its state-space, we present a new algorithm that computes deterministic upper and stochastic lower bounds of the value function similar to dual dynamic programming. We then show that the proposed algorithm terminates after a finite number of iterations. Finally, we demonstrate the efficacy of our approach on a high-dimensional numerical example from delivery slot pricing in attended home delivery.
\end{abstract}

\begin{keyword}
Dual dynamic programming; Function approximation; Real-time operations in transportation.
\end{keyword}

\end{frontmatter}
\thispagestyle{plain}
%===============================================================================
\section{Introduction}
Dynamic programming (DP) is an established tool to solve optimal control problems in a diverse set of applications. Since the exact DP approach typically scales poorly with the dimensionality of the problem at hand, various approximation schemes have been developed over the years. For the case of linear multistage stochastic systems, for example, stochastic dual DP provides a remedy to the ``curse of dimensionality'' by constructing successively tighter upper and lower bounds to the exact value function of the DP \citep{PEREIRA1991}. Moreover, this approach is able to provide stochastic guarantees on the performance of the approximation. However, current developments are mainly concerned with linear systems with piecewise-affine dynamics in \cite{PEREIRA1991} and in \cite{SHAPIRO2011} or piecewise-quadratic value functions in \cite{WARRINGTONETAL2019}. To the best of our knowledge, research on problems with discrete state-spaces has attracted only limited interest to date, see \cite{ZOU2019}.

In this paper, we extend the stochastic dual DP approach to value functions that are both submodular and concave extensible over a discrete domain. We define these properties more formally in Section \ref{sec:problem}. For any DP whose value function has these properties, we present a new algorithm that computes deterministic upper and stochastic lower bounds to the exact value function similar to stochastic dual DP theory.
	
Value functions with these properties arise, for example, in the revenue management problem in attended home delivery, i.e.\ the problem of finding the optimal pricing policy for a finite number of capacity-constrained delivery options. For this problem, \cite{LEBEDEVETAL2019B} show that the value function is concave extensible and submodular under certain conditions. We provide a numerical example with parameters that stems from this application to show the effectiveness of our new algorithm. Our approach advances the approximate DP approach of \citet{YANG2017}, who approximate the value function of the DP with an affine function, while our approximation takes the more expressive form of a piecewise affine function.
	
The paper is structured as follows:	In the remainder of this section, we introduce some notational conventions. Section \ref{sec:problem} formulates our problem of interest and the assumptions that our work builds upon. In Section \ref{sec:algorithm}, we present a novel algorithm to compute approximately optimal policies for discrete value functions under the imposed assumptions. Section \ref{sec:theoreticalresults} derives deterministic upper bounds and stochastic lower bounds to the exact value function and shows convergence of the algorithm in a finite number of iterations. In Section \ref{sec:gbdp_example}, we present a numerical example on a high-dimensional problem that cannot be solved by direct computation. Finally, Section \ref{sec:conclusions} concludes the paper and provides some directions for future research.

\emph{Notation:} Given some $s\in \mathbb{N}$, let $1_s$ be a column vector of all zeros apart from the $s$-th entry, which equals 1. Furthermore, let $1_0$ be a vector of zeros, i.e.\ $1_0:=\mathbf{0}$, and let $\mathbf{1}$ denote a vector of ones. Let $\innerprod{\cdot}{\cdot}$ denote the standard inner product of its arguments. Let $\Pr$ denote the probability operator and $\mathbb{E}$ the associated expectation.

\section{Problem statement}\label{sec:problem}
We consider a discrete-space, discrete-time, finite horizon DP. Define discrete states $x \in X \subset \mathbb{Z}^n$ and continuous and/or discrete decision variables $d\in D \subset \mathbb{Z}^a\times \Re^b$. Define the set $S:=\{1,2,...,n\}$. Define the transition probability between two states $x$ and $y$ under decision $d$ as $P_{x,y}(d)$, where we require $P_{x,y}(d)\geq 0$ for all $(x,y,d)\in X\times X\times D$. For all $x \in X$, we impose that $\sum_{y \in Y_+(x)}P_{x,y}(d)=1$, where $Y_+(x):=\{x+1_s\}_{s\in S \cup \{0\}}$. This requirement implies that transitions in $x$ are only possible in the positive direction and by at most a unit step along one dimension. Such models are typical for order-taking processes in revenue management (see \cite{ASDEMIRETAL2009}, \cite{SUH2011}, \cite{YANG2016} and \cite{YANG2017}). Furthermore, we define a finite time horizon \mbox{$T:=\{1,2,\dots,\bar{t}\}$}, a stage revenue function $g: \mathbb{Z}^{n}\times \mathbb{Z}^{n}\times (\mathbb{Z}^{a}\times\Re^b) \to \Re$ and a terminal cost function $C: \mathbb{Z}^n \to \Re$ to construct a DP of the following form:
\begin{align}\label{eq:dp}
V_t(x) :=&\; \underset{d\in D}{\max}\left\{\sum_{y\in Y_+(x)}P_{x,y}(d)\left(g(x,y,d)+V_{t+1}(y)\right)\right\}\nonumber\\ 
&\;\forall (x,t)\in X\times T,\text{ where } \nonumber \\
&\; V_{\bar{t}+1}(x)=-C(x) \quad \forall x\in X.
\end{align}
It is not strictly necessary for $g$ to be independent of $t$ as long as the assumptions stated below can be satisfied. However, as our interest lies in time-independent problems and to ease notation, we ignore time-dependency of $g$ in this paper. To represent the DP in a more compact form, we define the operator $\mathcal{T}$ through the relationship
\begin{equation}\label{eq:dpcompact}
V_t := \mathcal{T}V_{t+1}, \text{ for all } t\in T.
\end{equation}
We next introduce several definitions needed to express the assumptions that we impose on the DP in \eqref{eq:dp}. 
\begin{definition}\label{de:submod}
A function $f: \mathbb{Z}^n\to \Re$ is said to be \emph{submodular} if it satisfies
\begin{equation}\label{eq:submodular}
f(\max(y,z))+f(\min(y,z)) \leq f(y)+f(z)
\end{equation}
for all $(y,z) \in \mathbb{Z}^n \times \mathbb{Z}^n$, where the maximum and minimum are taken elementwise.
\end{definition}
 The following two definitions are commonly used in discrete convex analysis:
\begin{definition}\label{de:conc_clo}
Let $a \in \Re^n$ and $b \in \Re$. Then the \emph{concave closure} $\tilde{f}: \Re^n \to \Re \cup -\infty$ of a function \mbox{$f:\mathbb{Z}^n \to \Re \cup -\infty$} is defined as \cite[(2.1)]{MUROTA2001} $\tilde{f}(x):=\underset{a,b}{\inf}\set{a^\intercal x + b}{a^\intercal y + b \geq f(y) \hspace{4mm} \forall y \in \mathbb{Z}^n} \hspace{4mm}\forall x\in \text{dom}(\tilde{f}).$
\end{definition}
	
\begin{definition}\label{de:conc_ext}
A function $f: \mathbb{Z}^n \to \Re \cup -\infty$ is \emph{concave extensible} if and only if $f$ coincides with its concave closure $\tilde{f}$, i.e.\ $f(x)=\tilde{f}(x), \text{ for all } x\in \mathbb{Z}^n$ \citep[Lemma 2.3]{MUROTA2001}.
\end{definition}
	
These definitions allow us to state the assumptions that we impose on the DP in \eqref{eq:dp}:
\begin{assumption}\label{as:terminal}
The function $C$ is submodular and concave extensible in $x$.
\end{assumption}
\begin{assumption}\label{as:preserve}
The Bellman operator $\mathcal{T}$ preserves concave extensibility and submodularity of any concave extensible and submodular value function, i.e. if $V_{t+1}$ is submodular and concave extensible, then $V_t$ also has these properties for all $t\in T$.
\end{assumption}
In \citet[Theorem 2]{LEBEDEVETAL2019B}, it is shown that, under mild technical assumptions on the customer arrival rate, these assumptions are satisfied for the revenue management problem considered in Section \ref{sec:gbdp_example}.
	
\section{Proposed approximation algorithm}\label{sec:algorithm}
We first state our proposed approximation procedure in Algorithm \ref{alg:GBDP} below and subsequently describe the meaning of all symbols and the individual algorithm steps.

\begin{algorithm}[H]
\caption{Proposed approximation algorithm}\label{alg:GBDP}
\begin{algorithmic}[1]
\State Initialise parameters: $X, D, P_{x,y}, T, g, C$ and $i_{\max}$
\State Initialise $Q_t^{0}(x) \gets \infty$, for all $(x,t)\in X\times T$
\State Initialise $Q_{\bar{t}+1}^0(x) \gets -C(x)$, for all $x \in X$
\For{$i\in \{1,2,\dots, i_{\max}\}$}
	\State $x_1^i \gets \mathbf{0}$
	\For{$t\in T$} \Comment ``Forward sweep''
		\State \Longunderstack[l]{$d_t^i \gets d^* \in \underset{d \in D}{\text{argmax}}\left\{\sum_{x_{t+1}^i\in Y_+(x_t^i)}P_{x_t^i,x_{t+1}^i}(d) \right.$\\ \\$\left.\hphantom{d_t^i \gets}\vphantom{\sum_{x_t}^1Q^1}\times\left(g(x_t^i,x_{t+1}^i,d)+Q_{t+1}^{i-1}(x_{t+1}^i)\right)\right\}$}
		\State $x_{t+1}^i \gets x_{t}^i + \underset{x_{t+1}^i}{\text{sample}}\left\{P_{x_t^i,x_{t+1}^i}\left(d_t^i\right)\right\}$
	\EndFor
	\State $l(i)\gets \sum_{t=1}^{\bar{t}} g(x_t^i,x_{t+1}^i,d_t^{i}) - C(x_{\bar{t}+1})$
	\While{$t > 1$} \Comment ``Backward sweep''
		\State $Z(x_{t+1}^i) \gets \{x_{t+1}^i+1_s+1_s'\}_{(s,s')\in (S\cup \{\mathbf{0}\})\times (S\cup \{\mathbf{0}\})}$
		\If{$Q_{t+1}^{i-1}$ is submodular on $Z(x_{t+1}^i)$}
			\State $H^* \gets$ unique hyperplane through \hphantom{www} $\hphantom{wwwwww}\vphantom{sum^1}\left\{\left(y,(\mathcal{T}Q_{t+1}^{i-1})(y)\right)\right\}_{y \in Y_+(x_{t+1}^i)}$
		\Else
			\State $j^* \in \underset{j\in J_{t+1}^{i-1}}{\text{argmin}}\left\{\left(\mathcal{T}H_{t+1}^{j-1}\right)\left(x_{t+1}^i\right)\right\}$	
			\State $H^* \gets \mathcal{T}H_{t+1}^{j^*-1}$
  		\EndIf
  		\State $Q_t^{i}\gets \min\left\{H^*,Q_t^{i-1}\right\}$
  		\State $t \gets t-1 $ 
	\EndWhile
	\State $u(i)\gets Q_1^{i}(\mathbf{0})$ 
\EndFor
\end{algorithmic}
\end{algorithm}

Inspired by stochastic dual DP \citep{SHAPIRO2011}, the main idea of our algorithm is to alternate between generating sample paths in ``forward sweeps'' and refining the value function in ``backward sweeps''. The following sections describe this procedure in detail.
\subsection{Initialisation}
We first initialise all parameters of the DP in \eqref{eq:dp}. Denote the maximum number of iterations by $i_{\max} \in \mathbb{N}$ and let $I:=\{0,1,\dots,i_{\max}\}$. Let the value function approximation $Q_t^{i}$ for all $(i,t)\in I\times T$ be represented as the pointwise minimum of a finite number of affine functions, i.e.
\begin{equation}\label{eq:qminh}
Q_t^{i}(x):=\underset{j \in \{0,1,\dots,i\}}{\min}H_{t}^j(x), \text{ for all } x\in X,
\end{equation}
where $H_{t}^j: X \mapsto \Re$ describes a hyperplane, i.e.\ 
\begin{equation}
H_{t}^j(x):=\innerprod{a_{t}^j}{x} + b_{t}^j, \text{ for all } x\in X,
\end{equation}
where $a_{t}^j \in \Re^{n}, b_{t}^j \in \Re$ for all $(t,j)\in T\times I$. We characterise the set of supporting hyperplanes at $x$ as
\begin{equation}
J_t^{i}(x):=\underset{j \in \{0,1,\dots,i\}}{\text{argmin}}\left\{\innerprod{a_{t}^j}{x} + b_{t}^j \right\}
\end{equation}
for all $(x,i,t)\in X\times I \times T$. We construct $Q_t^i$ as a successively tighter upper bound of $V_t$ (as $i$ increases), i.e.\ $V_t(x)\leq Q_t^{i}(x)\leq Q_t^{i-1}(x)$ for all $(x,i,t)\in X\times I\times T$. In the $i$-th ''backward sweep``, $H_t^i$ is added to $Q_t^{i-1}$ for all $t\in T$ to form $Q_t^i$. To initialise $Q_t^0$, one could simply set $Q_t^0$ to be a single affine function with zero slope and infinite offset, which would guarantee that $Q_t^0$ is indeed an upper bound to $V_t$ for all $t\in T$. We discuss the possibility of closer initialisations in our numerical example in Section \ref{sec:gbdp_example}. We also initialise $Q_{\bar{t}+1}^i(x) := V_{\bar{t}+1}(x) = - C(x)$ for all $(x,i)\in X\times I$, which is a tight upper bound by the construction of the DP in \eqref{eq:dp}.
\subsection{``Forward sweep''}
Fix any iteration $i\in I\setminus\{0\}$. In each ``forward sweep'', we solve an approximate version of the Bellman equation in \eqref{eq:dp} forward in time, i.e.\ by replacing $V_t$ with its approximation from the previous iteration $Q_t^{i-1}$. By doing so, we compute suboptimal decisions $d_t^i$ for all $t \in T$ and simulate state transitions by sampling from the transition probability distribution given the  approximately optimal decisions. This defines a sample path $x_{t}^i $ for all $t \in T\cup \{\bar{t}+1\}$. At the end of each ``forward sweep'', we compute a stochastic lower bound on the total expected profit $V_1(\mathbf{0})$, which we denote by $l(i)$ for all $i \in I\setminus\{0\}$. We show that this is indeed a stochastic lower bound in Section \ref{sec:theoreticalresults}.
\subsection{``Backward sweep''}
Fix any iteration $i\in I\setminus\{0\}$. In each ``backward sweep'', we first check if $Q_{t+1}^{i-1}$ is submodular on $Z(x_{t+1}^i)$ by computing the sign of \eqref{eq:submodular} for all pairs of points $(y,y') \in Z(x_{t+1}^i) \times Z(x_{t+1}^i)$, such that $y\neq y'$. If the inequality in \eqref{eq:submodular} holds for all these points, we locally compute the exact DP stage problem on the set $Y_+(x_{t+1}^i)$ and then construct the hyperplane through $\left\{\left(y,(\mathcal{T}Q_{t+1}^{i-1})(y)\right)\right\}_{y \in Y_+(x_{t+1}^i)}$. Then, the resulting added hyperplane is necessarily an upper bound to $V_t(x)$ for all $x \in X$, which we show in Section \ref{sec:theoreticalresults}. In the opposite case, we need to compute a submodular upper bound on $Q_{t+1}^{i-1}$, which is readily given by the hyperplanes from which $Q_{t+1}^{i-1}$ is constructed. Therefore, we select the hyperplane $H_{t+1}^{j^*-1}$ that minimises the value at the evaluation point $x_{t}^{i}$, which therefore locally creates the tightest upper bound. It may be possible to construct other submodular upper bounds to non-submodular $Q_{t+1}^{i-1}$, however, steps 16 and 17 of Algorithm \ref{alg:GBDP} offer a simple implementation. Finally in step 18, we update the value function approximation as the pointwise minimum of the approximation from the previous iteration and the newly constructed hyperplane. We also compute an upper bound, $u(i)$ for all $i\in I\setminus\{0\}$, on the total expected profit $V_1(\mathbf{0})$. We show that this is indeed an upper bound in Section \ref{sec:theoreticalresults}.

\section{Theoretical results}\label{sec:theoreticalresults}
In this section, we show our main theoretical results on bounds on the exact value function and convergence properties of Algorithm \ref{alg:GBDP}.

\begin{prop}\label{pr:QgeqV}
Under Assumptions \ref{as:terminal} and \ref{as:preserve}, the approximate value function is an upper bound to the exact value function, i.e.\ $Q_t^i(x) \geq V_t(x)$ for all $(x,i,t)\in X\times I \times T$.
\end{prop}
\begin{pf}
We show this result by induction on $t$. In the base case, i.e.\ at the terminal condition, we have $Q_{\bar{t}+1}^i := V_{\bar{t}+1}(x) = - C(x)$ for all $i\in I$, which satisfies the proposition trivially by Assumption \ref{as:terminal}. Now assume for an induction hypothesis that  $Q_{t+1}^{i-1}(x) \geq V_{t+1}(x)$ for some $(i,t)\in I \times T$ and for all $x\in X$. Fix any $x$ in $X$ and distinguish the two cases of the if-statement in step 13 of Algorithm \ref{alg:GBDP}.
\subsubsection{Case I:}		
Suppose that $Q_{t+1}^{i-1}$ is submodular on $Z(x_{t+1}^i)$. Then $H^*$ is the hyperplane through the set of points $\{(y,(\mathcal{T}Q_{t+1}^{i})(y))\}_{y \in Y_+(x_{t+1}^i)}$. By \eqref{eq:qminh}, $Q_{t+1}^{i}$ is concave extensible since it is the pointwise minimum of a finite number of hyperplanes. Hence, we invoke Assumption \ref{as:preserve} to conclude that $\mathcal{T}Q_{t+1}^{i-1}$ is concave extensible and submodular. As shown by \citet[Appendix B.4]{LEBEDEVETAL2019B}, this implies that the hyperplane $H^*$ is a separating hyperplane, i.e.\ $H^*(x)\geq \mathcal{T}Q_{t+1}^{i-1}(x)$ for all $x\in X$. Define $d^V$ to be the maximiser of \eqref{eq:dp} and define $d^Q$ to be the maximiser of \eqref{eq:dp} with $V_{t+1}(y)$ replaced by $Q_{t+1}^{i-1}(y)$. We now show that the Bellman operator of the DP preserves the inequality $Q_{t+1}^{i-1}(x) \geq V_{t+1}(x)$, i.e.\ $\mathcal{T}Q_{t+1}^{i-1}(x) \geq \mathcal{T}V_{t+1}(x)$. To this end, fix $x\in X$ and consider				
\begin{align}
(\mathcal{T}Q_{t+1}^{i-1})(x)&= g(x,d^Q)+\sum_{y\in Y_+(x)}P_{x,y}(d^Q)Q_{t+1}^{i-1}(y)\nonumber\\
&\geq g(x,d^V)+\sum_{y\in Y_+(x)}P_{x,y}(d^V)Q_{t+1}^{i-1}(y)\nonumber\\
&\geq g(x,d^V)+\sum_{y\in Y_+(x)}P_{x,y}(d^V)V_{t+1}(y)\nonumber\\
&=(\mathcal{T}V_{t+1})(x),
\end{align}
where the first inequality follows from the supoptimality of $d^V$ for $(\mathcal{T}Q_{t+1}^{i-1})(x)$ and the second inequality follows from the induction hypothesis.
				
\subsubsection{Case II:}
Now suppose that $Q_{t+1}^{i-1}$ is not submodular on $Z(x_{t+1}^{i})$. Then
\begin{equation}\label{eq:hset}
H^* \in \set{\mathcal{T}H_{t+1}^{j-1}}{j\in J_{t+1}^{i-1}}.
\end{equation}
Furthermore, by \eqref{eq:qminh} and the induction hypothesis, we have 
\begin{equation}\label{eq:hqv}
H_{t+1}^{j-1}(x)\geq Q_{t+1}^{i-1}(x) \geq V_{t+1}(x), \text{ for all } (x,j)\in X\times J_{t+1}^{i-1}.
\end{equation}
We will now show that all elements of the set in \eqref{eq:hset} constitute upper bounds on $\mathcal{T}V_{t+1}$. To this end, fix any $(x,j)\in X\times J_{t+1}^{i-1}$. Define $d^H$ to be the maximiser of \eqref{eq:dp} with $V_{t+1}(y)$ replaced by $H_{t+1}^{j-1}(y)$. We can show that the Bellman operator of the DP preserves the inequality $Q_{t+1}^{i-1}(x) \geq V_{t+1}(x)$ using a similar argument as before, i.e.\
\begin{equation}
(\mathcal{T}H_{t+1}^{j-1})(x)\geq(\mathcal{T}V_{t+1})(x),
\end{equation}
which follows from the suboptimality of $d^V$ (see Case I) for $(\mathcal{T}H_{t+1}^{j-1})(x)$ and the fact that $H_{t+1}^{j-1}(x)\geq V_{t+1}(x)$ (see \eqref{eq:hqv}). Therefore, we conclude that $H^*(x)\geq \mathcal{T}V_t(x)$ for all $x\in X$ in the second case as well. 
		
Since both cases lead to an upper bound, i.e.\ $H^*(x)\geq \mathcal{T}V_{t+1}(x)$ for all $x\in X$, we infer that 
\begin{equation}
Q_{t}^i(x)= \min\left\{H^*(x),Q_{t+1}^{i-1}(x)\right\} \geq  \mathcal{T}V_{t+1}(x)
\end{equation}
for all $x\in X$. This concludes our induction argument and shows that $Q_t^i(x) \geq V_t(x)$ for all $(x,i,t)\in X\times I \times T$.\qed
\end{pf}

\begin{cor}\label{co:u}
Under Assumptions \ref{as:terminal} and \ref{as:preserve}, $u(i)$ is an upper bound to the total expected profit, i.e.\ $u(i) \geq V_1(\mathbf{0})$ for all $i\in I\setminus\{0\}$.
\end{cor}
\begin{pf}
This result follows immediately from Proposition \ref{pr:QgeqV} and by observing that $u(i)=Q_1^i(\mathbf{0})$ for all $i \in I\setminus\{0\}$ from step 22 of Algorithm \ref{alg:GBDP}.\qed
\end{pf}
\begin{prop}
In expectation, $l(i)$ is a lower bound to the total profit, i.e.\ $\mathbb{E}[l(i)]\leq V_1(\textbf{0})$, for all $i\in I\setminus\{0\}$.
\end{prop}
\begin{pf}
This must be the case since for any $i\in I\setminus\{0\}$, the value of $l(i)$ is obtained from suboptimal decisions $d_t^i$ for all $t\in T$, due to the use of  $Q_{t+1}^{i-1}$ instead of the exact (yet unavailable) $V_{t+1}$ in line 7 of Algorithm \ref{alg:GBDP}. It follows that $d_t^i$ is not a maximiser of the exact DP in \eqref{eq:dp}, which, by the principle of optimality, implies that the expected value accumulated under this suboptimal policy will not be greater than the value obtained under the optimal policy. Hence, $\mathbb{E}[l(i)]\leq V_1(\textbf{0})$ for all $i\in I\setminus \{0\}$. \qed
\end{pf}
%We can state a \textit{chance-constrained} lower bound on the profit generation efficiency of the algorithm for any single ``forward sweep'' under the following additional assumption.
%\begin{assumption}\label{as:sigma0}
%The probability of $\hat{\sigma}=0$ is zero.
%\end{assumption}
%Assumption \ref{as:sigma0} is rather weak as $l(k)$ for all $k\in K$ is highly unlikely to take identical values due to the typically high-dimensional state space and long time horizon. Should such degenerate cases still occur with significant probability, Assumption \ref{as:sigma0} can also be satisfied by increasing $k_{\max}$ until $l(k_{\max})$ takes a different value from $l(k)$ for all $k\in K\setminus \{k_{\max}\}$.
 
%\begin{prop}\label{pr:chance}
%Under Assumption \ref{as:sigma0}, $\eta_{\alpha}$ is a \textit{chance-constrained} lower bound on the profit generation efficiency of the algorithm, i.e.\ fix any $k\in K$, then $\eta_{\alpha}\leq \mathbb{E}[l(k)]/V_1(\mathbf{0})$ with probability $1-\alpha$.
%\end{prop}

The stochastic dual DP algorithm converges asymptotically in $i$ to the exact value function \citep{SHAPIRO2011}. We can strengthen this result for our algorithm by exploiting that the set of states $X$ is finite. Hence, the proposed algorithm converges in a finite number of steps under the following additional assumption.
\begin{assumption}\label{as:resampling}
If Algorithm \ref{alg:GBDP} produces a sample path in iteration $i$ that contains a state-time pair $(x,t)\in X\times T$, for which the exact value function has already been computed, i.e.\ if $Q_{t}^{i-1}(x_{t+1}^i)=V_{t}(x_{t+1}^i)$, then the state $x_{t+1}^{i}$ is resampled in step 8 by randomly selecting a state for which the exact value function has not yet been reached.
\end{assumption}
\begin{prop}\label{pr:converge}
Under Assumptions \ref{as:terminal}, \ref{as:preserve} and \ref{as:resampling}, the gap $u(i)-\mathbb{E}[l(i)]$ for all $i\in I\setminus\{0\}$ converges to $0$ in at most $\bar{t}|X|$ training iterations of Algorithm \ref{alg:GBDP}.
\end{prop}
\begin{pf}
We will show the proposition by induction on $t$. Consider the base case, when $Q_{\bar{t}+1}^0(x) = V_{\bar{t}+1}(x)$ for all $x\in X$. Then notice that in the ``backward sweep'', the proposed algorithm computes the Bellman equation from $\bar{t}+1 \to \bar{t}$ exactly for every $x\in X$. This is because $Q_{\bar{t}+1}^0$ is submodular by Assumption \ref{as:terminal} and hence, the if-statement in step 13 of Algorithm \ref{alg:GBDP} is true. By Assumption \ref{as:resampling}, $x_{\bar{t}+1}^i$ is resampled if for the time step transition $\bar{t}+1 \to \bar{t}$, the algorithm has not converged to the exact value function at $\bar{t}$ yet. Therefore, the value function is computed exactly at all $x\in X$ for the time step transition $\bar{t}+1 \to \bar{t}$ after at most $|X|$ iterations of the proposed algorithm, i.e.\ $Q_{\bar{t}}^{\hat{i}}(x) = V_{\bar{t}}(x)$ for all $x\in X$, where $\hat{i}\leq |X|$.

Now suppose by means of an induction hypothesis that for some $(t,i)\in T\times I$,  $Q_{t+1}^i(x) = V_{t+1}(x)$ for all $x\in X$. Then by Assumptions \ref{as:terminal} and \ref{as:preserve}, $V_{t+1}$ is submodular and hence, $Q_{t+1}^i$ is also submodular. By a similar argument to the base case, the proposed algorithm computes the exact value function for the time step transition $t+1 \to t$ in another $\hat{i}\leq |X|$ iterations. 

Hence, we conclude that for every time step transition, the proposed algorithm needs at most $|X|$ iterations to compute the exact value function for any one time step $t\in T$, which gives at most $\bar{t}|X|$ iterations for the total time horizon. Hence, after any $k\geq\bar{t}|X|$ iterations, $Q_{t}^{k}(x) = V_{t}(x)$ for all $(x,t)\in X\times T$. Therefore, both $\mathbb{E}[l(k)]=V_1(\mathbf{0})$, as well as $u(k)=Q_{1}^{k}(\mathbf{0})=V_1(\mathbf{0})$, which finally implies that  $\mathbb{E}[l(k)]=u(k)$ for all $k\geq\bar{t}|X|$ iterations.\qed
\end{pf}
Note that in practice it is likely to take an unacceptably large number of iterations for the algorithm to converge to the exact value function due to the large number of states $|X|$. Since the value function is computationally expensive to calculate for all states, we seek to generate closer approximations at points that are likely to be visited, i.e.\ points on the sample path, and to use this information to save on approximation accuracy for less likely samples.

Our ultimate objective is to solve problems with huge state spaces (e.g.\ $|X|\approx 10^{20}$) and long time horizons (e.g.\ $|T|\approx10^4$). In such scenarios, the need to resample the state as detailed in Assumption \ref{as:resampling} becomes negligible, because the required number of iterations to reach convergence is much larger than the maximum acceptable number of iterations. Therefore, from a practical point of view, we do not resample to satisfy Assumption \ref{as:resampling}. In this case, the proposed algorithm only asymptotically converges to the exact value function instead of in a finite number of steps, just as in stochastic dual DP \citep{SHAPIRO2011}.

\section{Numerical example}\label{sec:gbdp_example}
%\begin{figure*}[b]
%  \includegraphics[width=\textwidth,trim={0.6cm 0.6cm 0.6cm 0.6cm},clip]{PyBoxplots_02}
%  \caption{Boxplots of upper $u$ (blue/left) from 100 simulation runs and stochastic lower bounds $l(k)$ (orange/right) for one simulation run, but $k_{\max}=1000$ validation iterations. Iterations indicate the number of affine functions, from which the approximate value function is created.}
%  \label{fig:uplowbounds}
%\end{figure*}
We demonstrate our algorithm on a synthetic example of the revenue management in attended home delivery problem, where the objective is to price delivery time windows, called slots, dynamically over a finite time horizon to control the customer choice process in a profit-maximising way and such that all orders can still be fulfilled.

In this problem, $S$ is the set of delivery slots and the components of $x$ correspond to the number of orders placed in every delivery slot. The feasible set of states $X$ is defined by the maximum state vector $\bar{x}$, i.e.\ $X:=\set{x}{\mathbf{0}\leq x \leq \bar{x}}$. In the model, the prices $d_s$ for all $s\in S$ must be chosen from the interval $[\underline{d},\bar{d}]$ or $\infty$, which defines the feasible decision space $D:=\set{d}{d_s \in [\underline{d},\bar{d}]\cup\{\infty\}}$. Customers choose slots according to the multinomial logit model:
\begin{align}
P_{x,x}(d)&:=(1-\lambda)+\frac{\lambda}{\sum_{k\in S}\exp(\beta_c+\beta_k+\beta_dd_{k})+1},\nonumber\\
P_{x,x+1_s}(d)&:=\lambda \frac{\exp(\beta_c+\beta_s+\beta_d d_{s})}{\sum_{k\in S}\exp(\beta_c+\beta_k+\beta_dd_{k})+1}
\end{align}

for all $(x,d,s)\in X\times D \times S$, where $\lambda\in (0,1)$ is the probability that a customer arrives on the booking website, $\beta_c \in \Re$ denotes a constant offset, $\beta_s \in \Re$ represents a measure of the popularity for all delivery slots $s\in S$ and $\beta_d < 0$ is a parameter for the price sensitivity. This model is the same as in \citet{YANG2016}, who also detail how to estimate these parameters. The average revenue of an order is $r$ and the length of the time horizon, representing the booking period, is $\bar{t}$. The cost function $C$ represents the delivery cost for all lists of orders $x\in X$ accumulated at the end of the booking period. The challenge is to price the slots dynamically to maximise profits, which corresponds to solving a DP of the form of \eqref{eq:dp}, where $g(x,y,d):=r+d_s$ if $y=x+1_s$ for all $s \in S$ and otherwise, $g(x,y,d):=0$, i.e.\ the stage revenue is the average revenue plus delivery price for slot $s$ if slot $s$ is chosen and otherwise, it is zero. The DP in our numerical example takes the parameters in Table \ref{tab:alg01par} below, adapted from a real-world, multi-subarea case study by \citet{YANG2017} to a single delivery subarea scenario. Furthermore, we also adopt the customer choice parameters $\left(\beta_c,\beta_d,\{\beta_s\}_{s\in S}\right)$ from that paper.
\begin{table}[H]
\caption{Numerical example parameters.}
\begin{center}
\begin{tabular}{r|l}\label{tab:alg01par}
	$S$ & $\{1,2,\dots,17\}$ \\
	$\bar{x}$ & $[6,6,\dots,6]$ \\
	$\lambda$ & $0.008$ \\
	$\left[\underline{d},\bar{d}\,\right]$ & $[\textrm{\pounds}0,\textrm{\pounds}10]$ \\
	$r$ & \pounds$34.53$ \\
	$\bar{t}$ & $6990$ \\
	$C(x)$ & $\textrm{\pounds}0.083\times \mathbf{1}^{\intercal}x$ if $x\in X$ and $\infty$ otherwise
\end{tabular}
\end{center}
\end{table}
We have chosen $C(\mathbf{0})=0$, i.e.\ we ignore fixed cost, which has no implication on the computed pricing policy. Notice that to compute the value function directly, it would be necessary to compute the Bellman equation in \eqref{alg:GBDP} for all pairs $(x,t) \in X\times T$, which for the above parameters gives $(6+1)^{17}\times 6990 \approx 1.6\times 10^{18}$ evaluations. This is prohibitively large for any available computational technology. Therefore, an approximate algorithm is necessary.

For this type of DP, \citet[Theorem 2]{LEBEDEVETAL2019B} have shown that the Bellman operator preserves strict submodularity, i.e.\ the condition in \eqref{eq:submodular} holds with strict inequality, if a small enough $\lambda > 0$ is chosen. However, in this problem, the terminal condition $C$ is only weakly submodular in $x$. In fact, it is modular as it is an affine function of $x$. We assume that $\lambda=0.008$ from Table \ref{tab:alg01par} is small enough to satisfy Assumption \ref{as:preserve} in this problem set-up. We note however, that in the absence of strict submodularity, we compromise on the absolute theoretical guarantee that the upper bound of the proposed algorithm is indeed an upper bound to the exact value function.

To speed up computation, we initialise $Q_{t}^0$ for all $t\in T$ using the fixed point of the DP, $V^*$, which is a known upper bound to the exact value function computed at any $(x,t)\in X\times T$, i.e.\ $V^*(x)\geq V_t(x)$. This is always the case when $\mathcal{T}$ in \eqref{eq:dpcompact} is a monotone operator \cite[Chapter 3]{BERTSEKAS2012}. \cite{LEBEDEVETAL2019A} show that the fixed point is given analytically as
\begin{equation}\label{eq:fixedpoint}
V^*(x):= (\bar{d}+r)\mathbf{1}^\intercal (\bar{x}-x)-C(\bar{x})\text{ for all } x\in X.
\end{equation}
Hence, we use this result to set $Q_{t}^0(x)=V^*(x)$ instead of $\infty$ for all $(x,t)\in X\times T$. Note that the fixed point in \eqref{eq:fixedpoint} is an affine function, so the initialiser has low complexity, i.e.\ only one affine function describes $Q_{t}^0$.

We implement our algorithm in Julia \citep{BEZANSON2017} and run it for $i_{\max}=100$ iterations on an \mbox{i7-8565U} CPU at 1.80 GHz processor base frequency and with 16GB RAM, giving us a run time of 25 mins, 48 sec. In each iteration $i\in \{1,2,\dots i_{\max}\}$, we compute the upper bound on the expected profit $u(i)$ and the stochastic lower bound $l(i)$, corresponding to the sample profit obtained in a single ``forward sweep'' of Algorithm \ref{alg:GBDP}. We evaluate the cumulative moving average of the sample profits, i.e.\ $i^{-1}\sum_{j=1}^i l(j)$, which tends to the expected value of the stochastic lower bound as $i$ increases. Fig.\ \ref{fig:converge} shows how these bounds develop over 100 iterations. We make the following observations:

\begin{figure}[h]
\centering
  \includegraphics[width=0.42\textwidth,trim={0cm 0cm 0cm 0cm},clip]{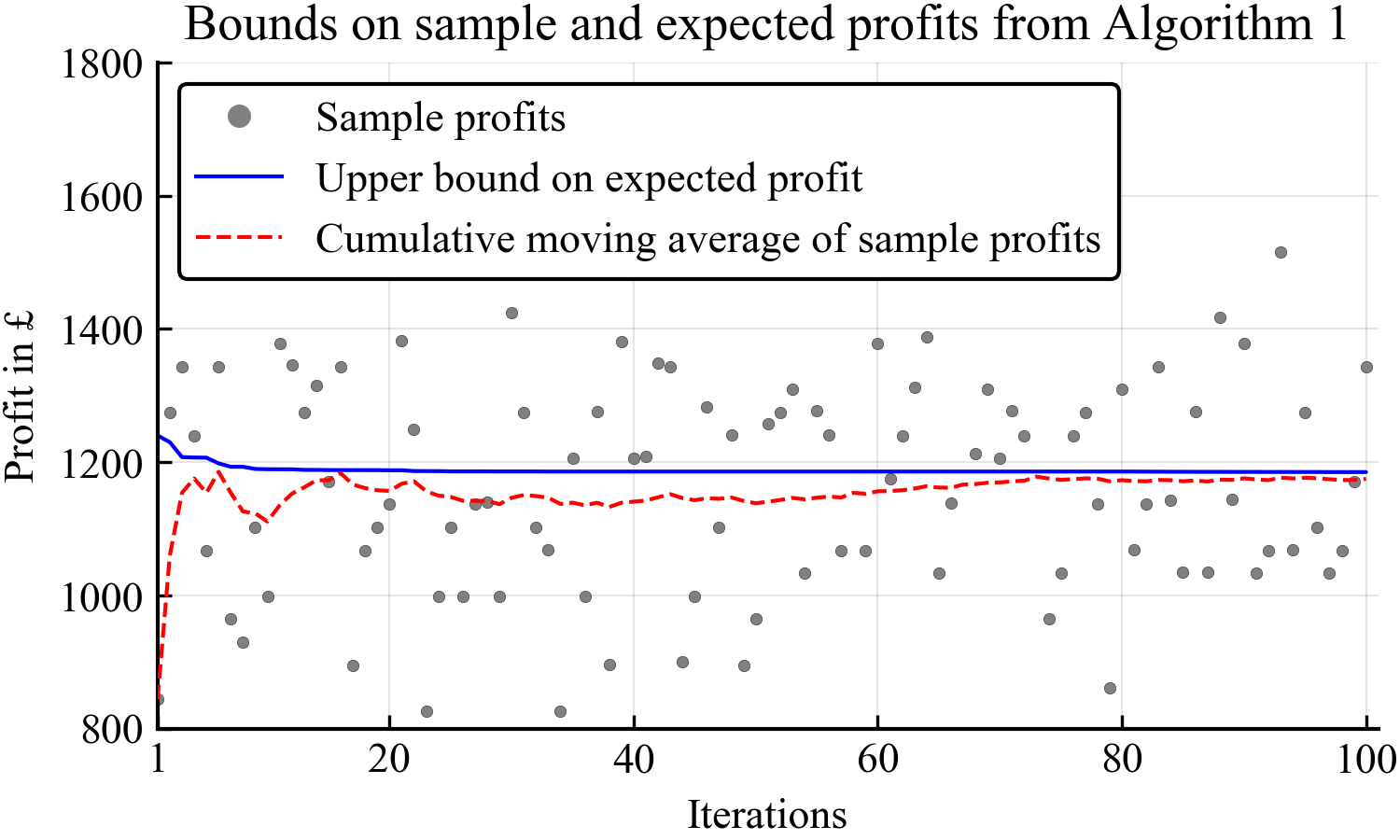}
  \caption{Plots of sample profits (grey dots), upper bound on expected profit $u(i)$ (blue solid line) and cumulative moving average of sample profits (red dashed line).}
  \label{fig:converge}
\end{figure}
%\begin{enumerate}

\emph{(1)} The upper bound converges within 10 iterations.

\emph{(2)} The cumulative moving average of the sample profits converges after about 100 iterations. One explanation for this is that Algorithm \ref{alg:GBDP} refines the value function approximation iteratively for all time steps. However, not all refinements propagate through all time steps of the DP. Therefore, each iteration step has a direct influence on the stochastic lower bound, which depends on the value function approximation at all time steps, while the upper bound might be unchanged as it only depends on the value function approximation at the first time step.

\emph{(3)} The sample profits have high variance at all algorithm iterations. This can be explained by the random customer arrival and choice process.

%\end{enumerate}
We investigate the influence of additional iterations, by comparing the performance of the pricing policies obtained after 1 and after 100 iterations of Algorithm \ref{alg:GBDP}. In particular, for both iteration counts, we simulate 1,000 booking periods (by using the ``forward sweep'' of Algorithm \ref{alg:GBDP}) and we compute the sample profits obtained in each period. The resulting histogram of sample profits is depicted in Fig. \ref{fig:histogram1} below, where we make the following observations:
\begin{figure}[h]
\centering
\subfloat[After 1 iteration.]{
  \includegraphics[width=0.42\textwidth,trim={0cm 0cm 0cm 0cm},clip]{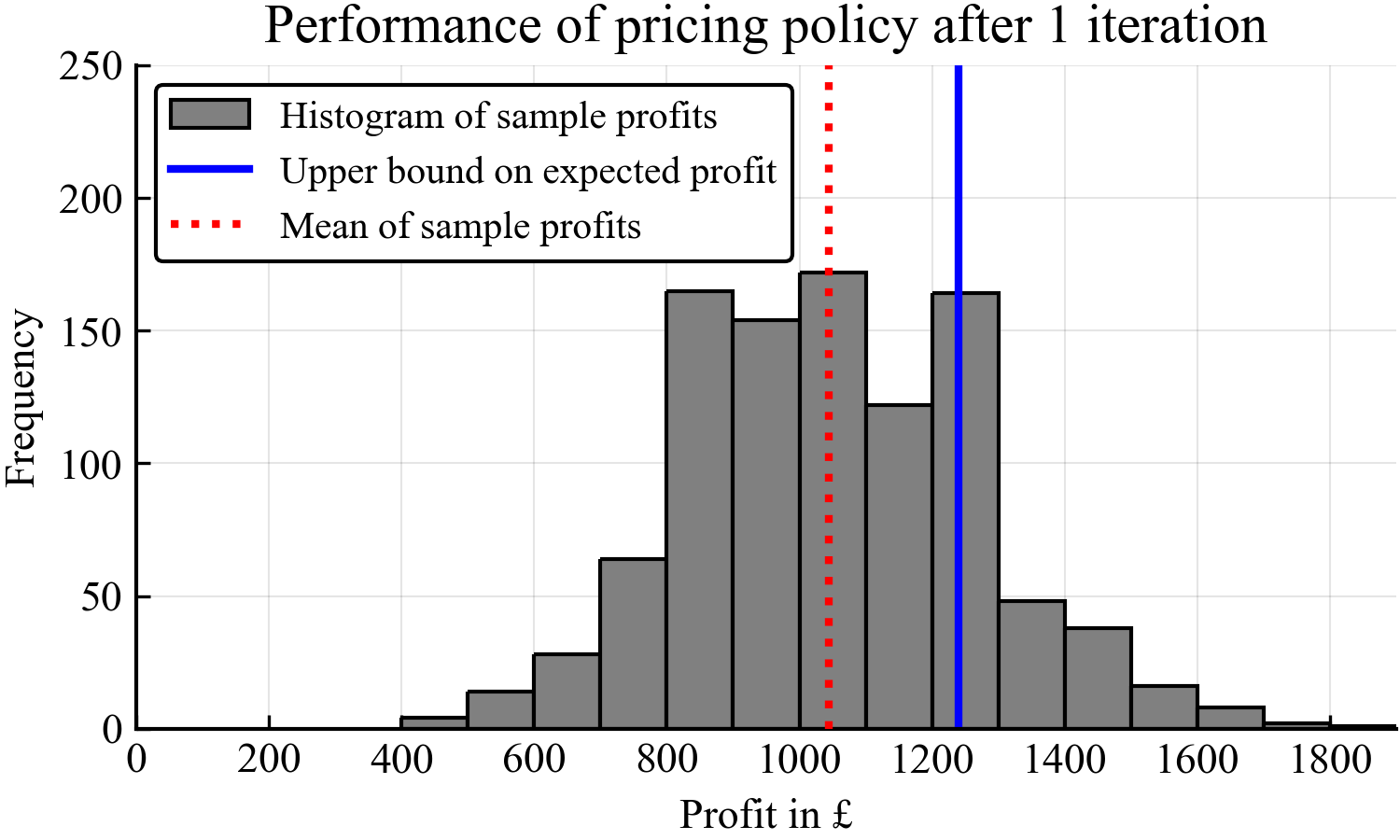}
 }
\hspace{0cm}
\subfloat[After 100 iterations.]{
  \includegraphics[width=0.42\textwidth,trim={0cm 0cm 0cm 0cm},clip]{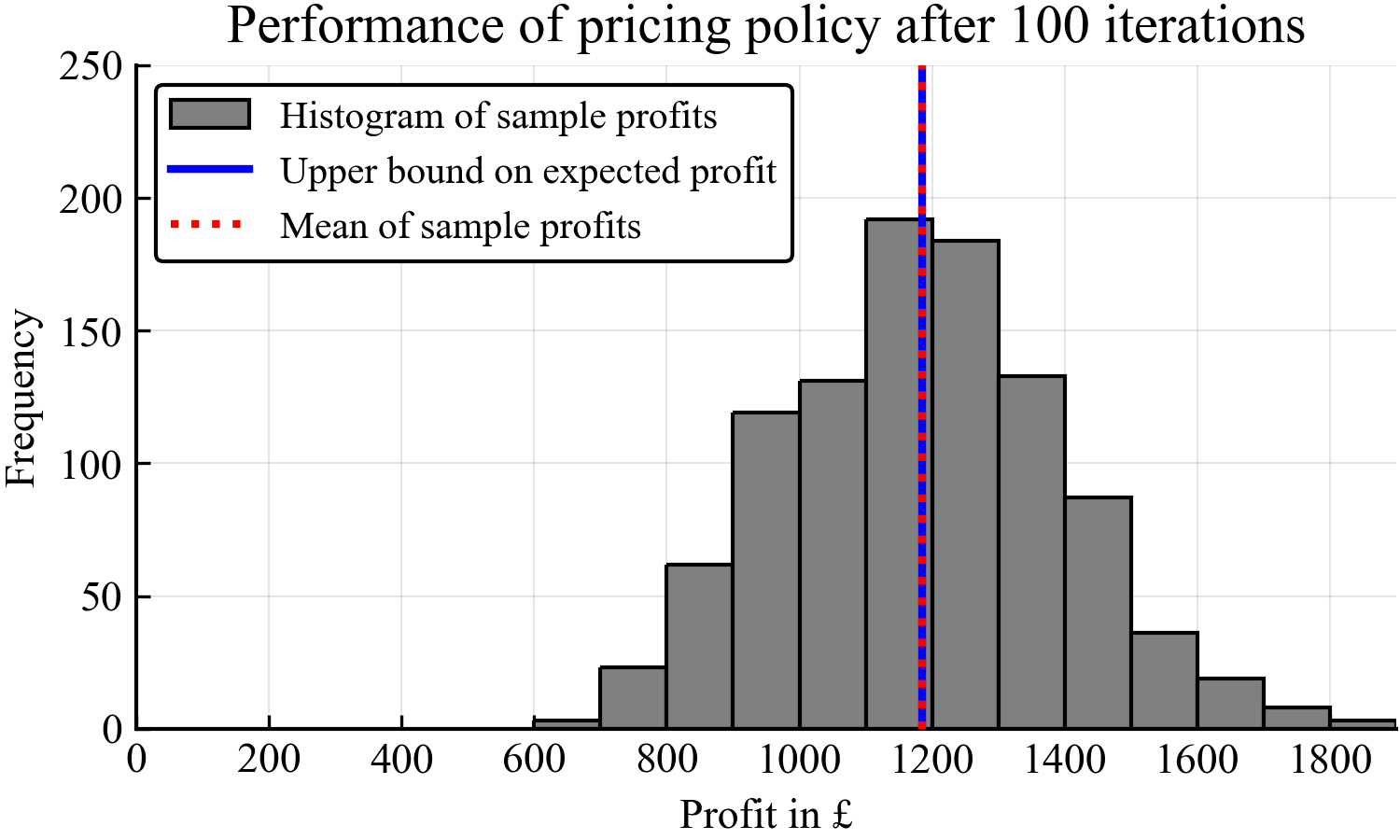}
 }
\caption{Histograms of 1,000 sample profits after 1 (a) and 100 (b) iterations of Algorithm \ref{alg:GBDP}.}
\label{fig:histogram1}
\end{figure}

\emph{(1)} The mean of the sample profits increases from \pounds 1,044 after 1 algorithm iteration to \pounds 1,185 after 100 iterations, an increase of 13.5\%.

\emph{(2)} The gap between upper bound on expected profit, $u(i)$, and the empirical mean of 1,000 samples of $l(i)$, decreases from \pounds 200 to approximately \pounds 0. 

\emph{(3)} The histograms show that the sample profits are more concentrated around their empirical mean after 100 iterations than after 1 iteration. This indicates that the variance of the sample profits can be decreased by increasing the number of algorithm iterations.

\section{Conclusions and future work}\label{sec:conclusions}
We present a new algorithm to compute approximate solutions to dynamic programs with submodular, concave extensible value functions. Similar to the results in stochastic dual dynamic programming, we derive deterministic upper and stochastic lower bounds to the exact value functions. Furthermore, we show that our algorithm converges in a finite number of iterations and we demonstrate our results by means of a numerical example of the revenue management problem in attended home delivery.
	
To the best of our knowledge, our new algorithm is the first to provide an upper bound on the expected profit for the class of problems under consideration. Comparing this upper bound with the stochastic lower bound, corresponding to the profit obtained in simulations of the dynamic program forward in time, we can quantify the profit generation efficiency of our algorithm. This quantity is a benchmark for other algorithms, possibly with weaker theoretical guarantees, but better performance in practice. 

Finally, the gap between upper bound and stochastic lower bound will allow to quantify and optimise the trade-off between quality of approximation and computational cost \textit{dynamically} as an application runs. For example, in the revenue management in attended home delivery problem, the terminal condition of the dynamic program is determined by the cost function, which is an approximation to the intractable capacitated vehicle routing problem with time windows (\cite{TOTH2014}). As time in the booking horizon progresses, orders come in, revealing the location and delivery times of customers. This aggregative information could be used to update the terminal condition of the DP, to re-run our proposed algorithm and hence, to update the pricing policy.

Directions for future work also include the derivation of probabilistic confidence intervals for the stochastic lower bound. Currently, the confidence intervals used in stochastic dual dynamic programming assume a Gaussian distribution for the sample profits, which holds only approximately for large sample sizes due to the central limit theorem \citep{SHAPIRO2011}. Non-parametric bounds, in the spirit of Chebyshev's inequality \citep{USPENSKY1937}, could help relax this assumption if they were adapted to cases, where the exact mean and variance of the underlying distribution are unknown and hence, are empirically estimated.

\begin{ack}
We gratefully acknowledge the helpful discussions with Michael Garstka, Department of Engineering Science, University of Oxford, on the Julia implementation of our Algorithm.
\end{ack}
\bibliography{Algorithm01_Ref02}

\end{document}